\documentclass[12pt]{article}
\usepackage[english]{babel}
\usepackage{epsfig}
\usepackage{amsfonts}
\usepackage{amssymb}
\usepackage{amsmath}
\usepackage{amsthm}
\usepackage{latexsym}
\usepackage{graphicx}
\usepackage{bm}
\usepackage{tabularx}
\usepackage{booktabs} 
\usepackage{enumerate}
\usepackage{caption}
\usepackage{wrapfig}  

\usepackage{mathtools}  
\usepackage{eqnarray}
\usepackage{enumitem}

\usepackage[titletoc]{appendix}

\usepackage[numbers]{natbib}
\usepackage{comment}
\usepackage{array} 
\usepackage{bbm} 
\usepackage{url}
\usepackage{subfig}
\usepackage{hyperref}
\usepackage{color}
\usepackage{float}

\usepackage{hyperref}
\hypersetup{colorlinks=true, citecolor=blue, linkcolor=red, urlcolor=green}
 
\newcolumntype{C}{>{$\displaystyle} c <{$}}
\usepackage{url}
\usepackage{color}
\usepackage{xcolor}
\usepackage[section]{placeins}

\theoremstyle{plain}
\newtheorem{theorem}{Theorem}[section]
\newtheorem{lemma}[theorem]{Lemma}
\newtheorem{corollary}[theorem]{Corollary}

\theoremstyle{definition}

\theoremstyle{remark}

\title{Overlap Times in the $GI^B/GI/\infty$ Queue}
\author{ Sergio Palomo   \\Systems Engineering \\ Cornell
University \\ {sdp85@cornell.edu} \and
Jamol Pender \footnote{Corresponding Author}
\\
School of Operations Research and Information Engineering
\\
Cornell University
\\
{jjp274@cornell.edu}
}

\usepackage{natbib}
\usepackage{graphicx}
\usepackage{amsmath}

\begin{document}

\maketitle

\begin{abstract}
Overlap times have been studied as a way of understanding the time of interaction between customers in a service facility.  Most of the previous analysis relies on the single jump assumption for arrivals, which implies the queue increases by one for each arrival epoch.  In this paper, we relax the single arrival assumption and explore the impact of having batch arrivals. Unfortunately, with batch arrivals it is not clear how one measures an overlap time between batches of customers.  Thus, we develop two ways of capturing the notion of an overlap time in a batch setting and derive exact results in the infinite server queue with batch arrivals.  Finally, we derive new results for analyzing overlap times of more than two batches.  
\end{abstract}

\section{Introduction}

The study of overlap times has emerged as a way of understanding how customers interact in service systems.  Recently, overlap times have been analyzed in infinite server queues, single server queues, and multi-server queues in \cite{kang2021queueing, palomooverlap, pender2021overlap, ko2022overlapping} respectively. This has been done mostly in Markovian cases since the analysis is a bit easier.  In \cite{kang2021queueing} the overlap times were used as a way of computing a new $R_0$ value for understanding infection rates in compartmentalized epidemic models in relation to COVID-19 transmission. \citet{palomooverlap} proves that the overlap distribution is exponential for the M/M/1 queue and \cite{pender2021overlap} studies the overlap distribution for infinite server queues that have renewal arrivals and i.i.d service times from a renewal process.  Finally, \citet{ko2022overlapping} uses fluid limits to approximate the overlapping time of customers in multi-server queues with time varying arrival rates.  

Beyond infectious disease transmission and design of service systems there are other applications where overlap times are important.  One application is public transportation like buses and subways.  In these transportation settings, several passengers get on the bus or subway in batches and ride for some time to their destination.  Thus, it is important to study the impact of batch arrivals to get a better understanding of how infectious disease spreads in these transportation settings.  

Fortunately, the study of batch arrival queues has increased significantly in recent years to due the emerging application of cloud-based data storage and processing. In this setting, the batches of jobs arriving to the system are collections of jobs submitted simultaneously by a user. These jobs are then served by each being processed individually and the results are returned to the user. For more discussion, see works such as \citet{lu2011join, pender2016law, xie2017pandas, yekkehkhany2018gb} and references therein.  Another relevant application in the context of infectious disease modeling is modeling the number infections of COVID-19.  See for example \citet{kaplan2020covid, morozova2020model, palomo2020flattening}.  In this setting, the results for patients who potentially have COVID-19 arrive in a large batch to be processed at a facility.  Moreover, the data that we observe from COVID-19 is also of batch form as counts are made daily.  Finally, an emerging application of batch queues is in context of autonomous vehicles moving in platoons (batches) down highways and roads, e.g.~\citet{mirzaeian2018queueing, hampshire2020beyond}. Such applications also serve as the inspiration for the batch arrival queue staffing problem studied by \citet{daw2021staff}.

In this work, we study overlap times for the $GI^B/GI/\infty$ queue.  Our results on the $GI^B/GI/\infty$ queue also serve as approximations and lower bounds for systems of finite capacity such as the $GI^B/GI/C$ queue. For example, consider the recent work on multi-server jobs, which are queueing systems where collections of jobs arrive together and also have a requirement that they must start together (see, e.g.,~\citet{rumyantsev2017stability, pender2017approximations, afanaseva2020stability, grosof2020stability, weng2020achieving, hong2021sharp, wang2021zero}, and references therein). The simultaneous start requirement is a salient model feature relative to batch arrival many server queues, but if there were infinitely many servers available then these models reduce to one another.  The multi-server setting model is known to be quite challenging to analyze, so we offer our following analysis in the unlimited server setting to understand the challenges and provide insights about more complicated models. 

In what follows we describe the contributions of our work and how the rest of the paper is organized.  

\subsection{Contributions of the Paper}
\begin{itemize}
    \item We derive explicit expressions for the steady state distribution for the overlap time for individual customers that are exactly $k$ batches apart.  
    \item We derive exact expressions for the tail distribution of the overlap time when we measure overlap in the first of each batch or last of each batch setting.  
    \item We also derive explicit expressions for the tail distribution of the overlap time when there are more than two customers or two batches.  
\end{itemize}

\subsection{Organization of Paper}
In Section \ref{Stoc:Model}, we describe the stochastic model that we will use in this work and the definition of overlap times.  We derive an exact expression for describing the overlap times for a pair of individual customers in the infinite server queue with batch arrivals i.e the $GI^B/GI/\infty$ queue.  We use this equation to compute the steady state distribution of the overlap time of customers that are $k$ batches apart.   In Section \ref{batch_perspective}, we also analyze the overlap times from a batch perspective. In particular we analyze two ways of comparing a batch.  The first is the first of each batch to leave and the second is the last of each batch to leave.  In Section \ref{more_than_two}, we analyze the overlap times of when there are more than two customers.  Finally in Section \ref{conclusion}, we provide a conclusion and some future research directions.  

\section{Infinite Server Overlap Times} \label{Stoc:Model}

In this section, we study the infinite server queue with the intention of understanding how much time do adjacent customers spend in the system together.  A similar type of analysis has been completed \citet{kang2021queueing, palomooverlap, pender2021overlap}.  Here we consider the $GI^B$/GI/$\infty$ queue where we let $A_i$ be the arrival time of the $i^{th}$ batch of customers, the inter-arrival time between the $i^{th}$ and $(i+1)^{th}$ batches is given by $A_{i+1} - A_i$.  We assume the inter-arrival times are i.i.d with cumulative distribution function (cdf) $H(x)$.  For convenience, we define $h_k(x)$ as the distribution of a sum of $k$ inter-arrival times.  We also assume that $S_{(i,j)}$ be the service time of the $i^{th}$ customer in the $j^{th}$ batch and the service times are i.i.d with cdf $G(x)$.   Thus, the departure time for the $j^{th}$ customer in the $n^{th}$ batch is given by the following expression

\begin{eqnarray}
D_{(j,n)} &=& S_{(j,n)} + A_{n}.
\end{eqnarray}
Given the arrival and departure time for each customer, it is possible to describe the overlap time between any pair of customers.  The overlap time between the $j^{th}$ customer in the $n^{th}$ batch and $\ell^{th}$ customer in the $(n+k)^{th}$ batch is given by
\begin{eqnarray}
O_{(j,n),(\ell,n+k)} &=& \left( \min( D_{(j,n)} , D_{(\ell,n+k)}  ) - A_{n+k} \right)^+ \\
&=& \left( \min( A_{n}  + S_{(j,n)}  , A_{n+k}  + S_{(\ell,n+k)}  ) - A_{n+k} \right)^+ \\
&=&  \left( (D_{(j,n)} - A_{n+k} )^+ \wedge S_{(\ell,n+k)}  \right) \\
&=&  \left( (S_{(j,n)} + A_{n} - A_{n+k} )^+ \wedge S_{(\ell,n+k)}  \right) \\
&=&  \left(  (S_{(j,n)} - \left( A_{n+k} - A_{n} \right) )^+ \wedge S_{(\ell,n+k)}  \right) .
\end{eqnarray}

\begin{theorem} \label{thm1}
Let $O_{j,k,\ell}$ have the steady state distribution of $O_{(j,n),(\ell,n+k)}$ in the $GI^B/GI/\infty$ queue, $\mathcal{S}$ and $\tilde{\mathcal{S}}$ be two independent service times with cdf $G(x)$, and let $\mathbb{P}(B = j) = p_j$,  then the tail distribution of $O_{j,k,\ell} = \lim_{n \to \infty} O_{(j,n),(\ell,n+k)}$ is given by

 \[
    \mathbb{P} \left( O_{j,k,\ell} > t  \right) = \left\{\begin{array}{lr}
       \overline{G}(t) , & \text{for } k=0,j = \ell, \ell \geq 1\\
        \overline{G}(t)^2  , & \text{for } k=0, j \neq \ell, \ell \geq 1\\
        \overline{G}(t)  \int^{\infty}_{0} \overline{G}\left( t + x \right) h_k(x) dx , & \text{for } k>0,j \neq \ell, \ell \geq 1
        \end{array}\right\}
  \]
where $h_k(x)$ is the density of the sum of $k$ i.i.d inter-arrival times.   
\begin{proof}
In the first condition where we have  $k=0,j = \ell, \ell \geq 1$, this corresponds to a self-overlap, which is precisely equal to the tail distribution of the service time.  The second condition where we have  $k=0,j \neq \ell, \ell \geq 1$, this corresponds to an overlap of two customers in the same batch.  Thus, the probability that they are both around by time $t$ is equal to the square of the tail distribution of the service time.  Finally, the last condition is exactly similar to the case with no batches, which is given in \citet{pender2021overlap}.  This completes the proof.
\end{proof}
\end{theorem}

Theorem \ref{thm1} provides the distribution of the case where we are comparing the overlap between \textbf{individual} customers.  However, we might be interested in the overlap of a batch itself.  In the sequel, we provide an analysis of this batch perspective.

\section{The Batch Perspective} \label{batch_perspective}

Unfortunately in the batch setting, there is not a unique way to define an overlap time between batches of customers.  In what follows, we describe two ways to define an overlap between customers in a batch.  Each one yields different expressions for calculating the tail distribution of the respective overlap time definition. We recognize that there are many more ways to define an interaction or overlap between a batch of customers, however, these two are quite natural.  However, before we derive the overlap times for the batch queues, we prove a simple result about the tail distribution of the maximum and minimum of a random sequence of random variables.  

\begin{lemma}
Let $B$ be a discrete random variable on the positive integers i.e. $\{ 1,2,3...\}$ with probabilities \{$p_1, p_2, ... p_n$\} and let $X = (X_1, X_2, ..., X_B)$ be a random vector with random size given by the random variable $B$.  Then, we have the following distributions for $\min(X)$ and $\max(X)$

\begin{eqnarray}
\mathbb{P} \left( \min(X) > t \right) &=& \sum^{\infty}_{j=1} p_j  \overline{G}(t)^j
\end{eqnarray}
and 
\begin{eqnarray}
\mathbb{P} \left( \max(X) > t \right) &=& 1 - \sum^{\infty}_{j=1} p_j  G(t)^j.
\end{eqnarray}
\begin{proof}
For the minimum we have that
\begin{eqnarray*}
\mathbb{P} \left( \min(X) > t \right) &=& \sum^{\infty}_{j=1} \mathbb{P} \left( \min(X) > t | B = j\right) \cdot  \mathbb{P} \left( B = j\right)  \\
&=& \sum^{\infty}_{j=1} \prod^{j}_{i=1} \mathbb{P} \left( X_i > t \right) \cdot  p_j  \\
% &=& \sum^{\infty}_{j=1}  \overline{G}(t)^j \cdot  p_j  \\
&=& \sum^{\infty}_{j=1} p_j  \overline{G}(t)^j.
\end{eqnarray*}
Finally, for the maximum we have that
\begin{eqnarray*}
\mathbb{P} \left( \max(X) > t \right) &=& 1 - \mathbb{P} \left( \max(X) \leq t \right) \\ 
&=& 1- \sum^{\infty}_{j=1} \mathbb{P} \left( \max(X) \leq t | B = j\right) \cdot  \mathbb{P} \left( B = j\right)  \\
&=& \sum^{\infty}_{j=1} \prod^{j}_{i=1} \mathbb{P} \left( X_i \leq t \right) \cdot  p_j  \\
% &=& 1 - \sum^{\infty}_{j=1}  G(t)^j \cdot  p_j  \\
&=& 1 - \sum^{\infty}_{j=1} p_j  G(t)^j.
\end{eqnarray*}
\end{proof}
\end{lemma}

The reader should note two things.  First, when the batch size is one with probability one, then this reduces to the directly to the tail distribution of the random variable and this holds for either the max or min since they are equivalent for one random variable.  Second, we note that when the random variable $B$ has finite support, then the infinite sums reduce to finite ones.  The finite support of the random variable $B$ is important from a computational perspective since it often can be the difference between an exact expression and an approximation of it. 

\subsection{Last of Each Batch to Leave}

The first method for capturing overlap between customers in batches is called "last of each batch to leave". This method defines an overlap for a batch as any time that any of customers of each batch overlap.  Mathematically, if we define $\vec{D}_n^{B_n} = ( D_{n,1}, D_{n,2}, ... , D_{n, B_n})$ as the vector of departure times of the $n^{th}$ batch and $\vec{S}_n^{B_n} = ( S_{n,1}, S_{n,2}, ... , S_{n, B_n})$ as the vector of service times of the $n^{th}$ batch, then we have the following expression of the overlap time of the "last of each batch to leave" 

\begin{eqnarray}
O_{n,n+k} &=& \left( \min \left( \max \left(\vec{D}_n^{B_n} \right) , \max \left(\vec{D}_{n+k}^{B_{n+k}} \right) \right) - A_{n+k} \right)^+ \\
&=&  \left(  \left(\max \left(\vec{S}_n^{B_n} \right) - \left( A_{n+k} - A_{n} \right) \right)^+ \wedge \max \left( \vec{S}_{n+k}^{B_{n+k}} \right) \right) .
\end{eqnarray}

\begin{theorem} \label{ssdist_max_gen}
Let $O_k$ be the distribution of $O_{n,n+k}$ for all values of $n$ in the $GI^B/GI/\infty$ queue and let $\mathcal{A}_k$, be the sum of k i.i.d inter-arrival times from a renewal process with density $h_k(x)$, then the tail distribution of $O_k$ under the last of each batch to leave overlap process is given by
\begin{eqnarray} \label{overlap_dist_det_last}
\mathbb{P} \left( O_{k} > t  \right) &=&    \left( 1 - \mathbb{E} \left[ G  \left( t + \mathcal{A}_k  \right)^\mathcal{B}  \right] \right) \cdot \left(  1 - \mathbb{E} \left[ G(t)^\mathcal{B}  \right]  \right).
\end{eqnarray}
\begin{proof}
First, we need to decompose the overlap probability into two probabilities by using a property of the minimum of two independent random variables i.e.
\begin{eqnarray*}
\mathbb{P} \left( O_{k} > t  \right) &=& \mathbb{P} \left(  \left(  \left(\max \left(\vec{S}_n^{B_n} \right) - \left( A_{n+k} - A_{n} \right) \right)^+ \wedge \max \left( \vec{S}_{n+k}^{B_{n+k}} \right) \right)   > t  \right) \\
&=& \mathbb{P} \left( \left(\max \left(\vec{S}_n^{B_n} \right) - \left( A_{n+k} - A_{n} \right) \right)^+ > t  \right) \cdot  \mathbb{P} \left(  \max \left( \vec{S}_{n+k}^{B_{n+k}} \right)   > t  \right) \\
&=& \mathbb{P} \left( \max \left(\vec{S}_n^{B_n} \right)   > t +  A_{n+k} - A_{n}  \right) \cdot  \mathbb{P} \left(  \max \left( \vec{S}_{n+k}^{B_{n+k}} \right)   > t  \right) \\
&=& \mathbb{E} \left[ 1 -  \sum^{\infty}_{j=1} p_j  G  \left( t +  A_{n+k} - A_{n}  \right)^j \right] \cdot \left(  1 - \sum^{\infty}_{j=1} p_j G \left(t \right)^j \right) \\
&=&  \left( 1 - \sum^{\infty}_{j=1} p_j \mathbb{E} \left[  G  \left( t +  A_{n+k} - A_{n}  \right)^j \right] \right) \cdot \left(  1 - \sum^{\infty}_{j=1} p_j G \left(t \right)^j \right) \\
% &=&  \sum^{\infty}_{j=1} p_j \mathbb{E} \left[  G  \left( t +  \mathcal{A}_k  \right)^j \right] \cdot \left(  1 - \sum^{\infty}_{j=1} p_j G \left(t \right)^j \right) \\
&=&  \left( 1 - \sum^{\infty}_{j=1} p_j \left( \int^{\infty}_{0}  G  \left( t +  x  \right)^j h_k(x) dx \right) \right) \cdot \left(  1 - \sum^{\infty}_{j=1} p_j G \left(t \right)^j \right) \\
&=&  \left( 1 - \mathbb{E} \left[ G  \left( t + \mathcal{A}_k  \right)^\mathcal{B}  \right] \right) \cdot \left(  1 - \mathbb{E} \left[ G(t)^\mathcal{B}  \right]  \right).
\end{eqnarray*}
This completes the proof.
\end{proof}
\end{theorem}

\begin{corollary} \label{ssdist_max}
Let $O_k$ be the distribution of $O_{n,n+k}$ for all values of $n$  in the $M^\mathcal{B}/M/\infty$ queue, then the tail distribution of $O_k$ under the last of each batch to leave setting is given by
\begin{eqnarray} \label{overlap_dist_last}
\mathbb{P} \left( O_{k} > t  \right) &=&  \left( 1 - \mathbb{E} \left[\sum^{\mathcal{B}}_{m=0} {\mathcal{B} \choose m} (-1)^m e^{-\mu t m}  \left( \frac{\lambda}{\lambda + \mu m } \right)^k \right] \right) \cdot \left( 1 - \mathbb{E} \left[ ( 1- e^{-\mu t  })^\mathcal{B}  \right]\right) .  
\end{eqnarray}
\begin{proof}
Using the previous result in Theorem \ref{ssdist_max_gen}, we have that
\begin{eqnarray*}
\lefteqn{ \mathbb{P} \left( O_{k} > t  \right) } \\ &=& \left( 1 - \sum^{\infty}_{j=1}  p_j \int^{\infty}_{0} \left( 1 - e^{- \mu (t+x) } \right)^j \frac{ (\lambda x)^{k-1}}{(k-1)!} \lambda e^{-\lambda x} dx \right) \cdot \left( 1 - \sum^{\infty}_{j=1}  p_j ( 1- e^{-\mu t})^j \right)\\
&=& \left( 1 - \sum^{\infty}_{j=1}  p_j \int^{\infty}_{0} \left( \sum^{j}_{m=0} {j \choose m} (-1)^{m} e^{-\mu (t +x) m}  \right) \frac{ (\lambda x)^{k-1}}{(k-1)!} \lambda e^{-\lambda x} dx \right) \cdot \left( 1 - \sum^{\infty}_{j=1}  p_j ( 1- e^{-\mu t})^j \right)\\
&=& \left( 1 - \sum^{\infty}_{j=1}  p_j \sum^{j}_{m=0} {j \choose m} (-1)^{m} \int^{\infty}_{0}  e^{-\mu (t +x) m}  \frac{ (\lambda x)^{k-1}}{(k-1)!} \lambda e^{-\lambda x} dx \right) \cdot \left( 1 - \sum^{\infty}_{j=1}  p_j ( 1- e^{-\mu t})^j \right)\\
&=& \left( 1 - \sum^{\infty}_{j=1}  p_j \sum^{j}_{m=0} {j \choose m} (-1)^m e^{-\mu t m}  \left( \frac{\lambda}{\lambda + \mu m } \right)^k \right) \cdot \left( 1 - \sum^{\infty}_{j=1}  p_j ( 1- e^{-\mu t})^j \right)\\
% &=&  \left( 1 - \mathbb{E} \left[\sum^{\mathcal{B}}_{m=0} {\mathcal{B} \choose m} (-1)^m e^{-\mu t m}  \left( \frac{\lambda}{\lambda + \mu m } \right)^k \right] \right) \cdot \left( 1 - \sum^{\infty}_{j=1}  p_j ( 1- e^{-\mu t})^j \right)\\
&=&  \left( 1 - \mathbb{E} \left[\sum^{\mathcal{B}}_{m=0} {\mathcal{B} \choose m} (-1)^m e^{-\mu t m}  \left( \frac{\lambda}{\lambda + \mu m } \right)^k \right] \right) \cdot \left( 1 - \mathbb{E} \left[ ( 1- e^{-\mu t  })^\mathcal{B}  \right]\right) .  
\end{eqnarray*}
This completes the proof.
\end{proof}
\end{corollary}

\begin{corollary} \label{ssdist_max_det}
Let $O_k$ be the distribution of $O_{n,n+k}$ for all values of $n$  in the $GI^B/D/\infty$ queue, then the tail distribution of $O_k$ under the last of each batch to leave setting is given by
\begin{eqnarray} \label{overlap_dist_det}
\mathbb{P} \left( O_{k} > t  \right) &=&   \mathbb{F}^{(k)} \left( (\Delta - t)^+  \right)
\end{eqnarray}
where $\mathbb{F}^{(k)}(\cdot)$ is the $k$-fold convolution of the inter-arrival distribution $F(x)$.
\begin{proof}
Using the previous result in Theorem \ref{ssdist_max_gen}, we have that
\begin{eqnarray*}
\mathbb{P} \left( O_{k} > t  \right) &=& \mathbb{P} \left(  \left(  \left(\max \left(\vec{S}_n^{B_n} \right) - \left( A_{n+k} - A_{n} \right) \right)^+ \wedge \max \left( \vec{S}_{n+k}^{B_{n+k}} \right) \right)   > t  \right) \\
&=& \mathbb{P} \left( \left( \Delta - \left( A_{n+k} - A_{n} \right) \right)^+ > t  \right) \cdot  \mathbb{P} \left(  \Delta   > t  \right) \\
&=& \mathbb{F}^{(k)} \left( (\Delta - t)^+  \right) .
\end{eqnarray*}
This completes the proof.
\end{proof}
\end{corollary}

\subsection{First of Each Batch to Leave}

The second method for capturing overlap between customers in batches is called "first of each batch to leave".  This method defines an overlap for a batch as any time that \textbf{all} of the customers of each batch overlap.  Mathematically, we have the following expression of the overlap time of the "first of each batch to leave" 

\begin{eqnarray}
O_{n,n+k} &=& \left( \min \left( \min \left(\vec{D}_n^{B_n} \right) , \min \left(\vec{D}_{n+k}^{B_{n+k}} \right) \right) - A_{n+k} \right)^+ \\
&=&  \left(  \left(\min \left(\vec{S}_n^{B_n} \right) - \left( A_{n+k} - A_{n} \right) \right)^+ \wedge \min \left( \vec{S}_{n+k}^{B_{n+k}} \right) \right)  .
\end{eqnarray}

\begin{theorem} \label{ssdist_min_gen}
Let $O_k$ be the distribution of $O_{n,n+k}$ for all values of $n$ in the $GI^B/GI/\infty$ queue and $\mathcal{A}_k$, be the sum of k i.i.d inter-arrival times from a renewal process with density $h_k(x)$, then the tail distribution of $O_k$ under the first of each batch to leave setting is given by
\begin{eqnarray} \label{overlap_dist_det_first}
\mathbb{P} \left( O_{k} > t  \right) &=&   \mathbb{E} \left[ \overline{G}(t+\mathcal{A}_k)^\mathcal{B}  \right]\cdot  \mathbb{E} \left[ \overline{G}(t)^\mathcal{B}  \right] .
% \sum^{\infty}_{j=1} p_j \left( \int^{\infty}_{0}  \overline{G}  \left( t +  x  \right)^j h_k(x) dx \right) \cdot \left(   \sum^{\infty}_{j=1} p_j \overline{G} \left(t \right)^j \right) .
\end{eqnarray}
\begin{proof}
First, we need to decompose the overlap probability into two probabilities by using a property of the minimum of two independent random variables i.e.
\begin{eqnarray*}
\mathbb{P} \left( O_{k} > t  \right) &=& \mathbb{P} \left(  \left(  \left(\min \left(\vec{S}_n^{B_n} \right) - \left( A_{n+k} - A_{n} \right) \right)^+ \wedge \min \left( \vec{S}_{n+k}^{B_{n+k}} \right) \right)   > t  \right) \\
&=& \mathbb{P} \left( \left(\min \left(\vec{S}_n^{B_n} \right) - \left( A_{n+k} - A_{n} \right) \right)^+ > t  \right) \cdot  \mathbb{P} \left(  \min \left( \vec{S}_{n+k}^{B_{n+k}} \right)   > t  \right) \\
&=& \mathbb{P} \left( \min \left(\vec{S}_n^{B_n} \right)   > t +  A_{n+k} - A_{n}  \right) \cdot  \mathbb{P} \left(  \min \left( \vec{S}_{n+k}^{B_{n+k}} \right)   > t  \right) \\
&=& \mathbb{E} \left[  \sum^{\infty}_{j=1} p_j  \overline{G}  \left( t +  A_{n+k} - A_{n}  \right)^j \right] \cdot \left(   \sum^{\infty}_{j=1} p_j \overline{G} \left(t \right)^j \right) \\
&=&  \sum^{\infty}_{j=1} p_j \mathbb{E} \left[  \overline{G}  \left( t +  A_{n+k} - A_{n}  \right)^j \right] \cdot \left(   \sum^{\infty}_{j=1} p_j \overline{G} \left(t \right)^j \right) \\
&=&  \sum^{\infty}_{j=1} p_j \mathbb{E} \left[  \overline{G}  \left( t +  \mathcal{A}_k  \right)^j \right] \cdot \left(   \sum^{\infty}_{j=1} p_j \overline{G} \left(t \right)^j \right) \\
&=&  \sum^{\infty}_{j=1} p_j \left( \int^{\infty}_{0}  \overline{G}  \left( t +  x  \right)^j h_k(x) dx \right) \cdot \left(   \sum^{\infty}_{j=1} p_j \overline{G} \left(t \right)^j \right)\\
&=& \mathbb{E} \left[ \overline{G}(t+\mathcal{A}_k)^\mathcal{B}  \right]\cdot  \mathbb{E} \left[ \overline{G}(t)^\mathcal{B}  \right].
\end{eqnarray*}
This completes the proof.
\end{proof}
\end{theorem}

\begin{corollary} \label{ssdist_min}
Let $O_k$ be the distribution of $O_{n,n+k}$ for all values of $n$ in the $M^B/M/\infty$ queue and $\mathcal{A}_k$, be the sum of k i.i.d inter-arrival times from a renewal process with density $h_k(x)$, then the tail distribution of $O_k$ under the first of each batch to leave setting is given by
\begin{eqnarray} \label{overlap_dist_first}
\mathbb{P} \left( O_{k} > t  \right) &=&   \mathbb{E} \left[ e^{- \mu t \mathcal{B} } \left( \frac{\lambda}{\lambda + \mu \mathcal{B}} \right)^k \right] \cdot  \mathbb{E} \left[ e^{- \mu t \mathcal{B} }  \right]  .
\end{eqnarray}
\begin{proof}
Using the previous result in Theorem \ref{ssdist_min_gen}, we have that
\begin{eqnarray*}
\mathbb{P} \left( O_{k} > t  \right) 
% &=& \mathbb{P} \left(  \left(  \left(\min \left(\vec{S}_n^{B_n} \right) - \left( A_{n+k} - A_{n} \right) \right)^+ \wedge \min \left( \vec{S}_{n+k}^{B_{n+k}} \right) \right)  > t  \right) \\
% &=& \mathbb{P} \left( \left(\min \left(\vec{S}_n^{B_n} \right) - \left( A_{n+k} - A_{n} \right) \right)^+ > t  \right) \cdot  \mathbb{P} \left(  \min \left( \vec{S}_{n+k}^{B_{n+k}} \right)   > t  \right) \\
% &=& \mathbb{P} \left( \min \left(\vec{S}_n^{B_n} \right)   > t +  A_{n+k} - A_{n}  \right) \cdot  \mathbb{P} \left(  \min \left( \vec{S}_{n+k}^{B_{n+k}} \right)   > t  \right) \\
&=& \mathbb{E} \left[  \sum^{\infty}_{j=1} p_j e^{- \mu j ( t +  A_{n+k} - A_{n}) }  \right] \cdot \left(  \sum^{\infty}_{j=1} p_j  e^{-\mu j t} \right) \\
&=&   \sum^{\infty}_{j=1} p_j e^{- \mu j  t} \mathbb{E} \left[ e^{- \mu j (A_{n+k} - A_{n}) }  \right] \cdot \left(  \sum^{\infty}_{j=1} p_j  e^{-\mu j t} \right) \\
&=&   \sum^{\infty}_{j=1} p_j e^{- \mu j  t} \mathbb{E} \left[ e^{- \mu j \mathcal{A} }  \right]^{k} \cdot \left(  \sum^{\infty}_{j=1} p_j  e^{-\mu j t} \right) \\
% &=&  \sum^{\infty}_{i=1} \sum^{\infty}_{j=1} p_j p_i  e^{- \mu (i+j)  t} \mathbb{E} \left[ e^{- \mu j \mathcal{A} }  \right]^{k} \\
&=&  \mathbb{E} \left[ e^{- \mu t \mathcal{B} } \left( \frac{\lambda}{\lambda + \mu \mathcal{B}} \right)^k \right] \cdot  \mathbb{E} \left[ e^{- \mu t \mathcal{B} }  \right] .
\end{eqnarray*}
This completes the proof.
\end{proof}
\end{corollary}

\section{More Than A Pair Setting} \label{more_than_two}

Our above analysis for overlap times only characterizes the overlap for two customers or two batches. In some settings it may be interesting to characterize the behavior of more than a pair of customers or a pair of batches.  We analyze the setting of more than two customers in what follows below.  

\subsection{Individual Customer Perspective}

For our first result, we characterize the overlap time of an arbitrary number of customers in the $GI/GI/\infty$ queue.  In particular, the overlap time for a k-tuple of customers in the $GI/GI/\infty$ queue, which is denoted by $O_{n_1, n_2, \dots , n_k}$, can be computed by the following expressions 

\begin{eqnarray*}
O_{n_1, n_2, \dots , n_k} &=& \left( \min( D_{n_1} , D_{n_2}, \cdots , D_{n_k} ) - A_{n_k} \right)^+ \\
&=& \left( \min( A_{n_1} + S_{n_1} , \dots , A_{n_j} + S_{n_j},  A_{n_k} + S_{n_k} ) - A_{n_k} \right)^+ \\
&=& (S_{n_1} + A_{n_1} - A_{n_k} )^+ \boldsymbol{\wedge} \cdots \boldsymbol{\wedge} (S_{n_j} + A_{n_j} - A_{n_k} )^+ \boldsymbol{\wedge} S_{n_k}  .
\end{eqnarray*}

First, note that the distance between any two pairs of customers is arbitrary.  Moreover, unlike the two customer case, we observe that the random variables are not independent anymore.  This is because the inter-arrival times are connected in each of the customers except for the last one, which only has a service time.  Nonetheless, in the exponential setting, we can still show that the arrival process only serves to determine the constants and not the functional form of the overlap distribution. 

\begin{theorem} \label{tuple_result}
Let $O_{n_1, n_2, \dots , n_k}$ be the distribution of the overlap time for customers within distance $n_k - n_1$ in the $GI/GI/\infty$ queue, then the tail distribution of $O_{n_1, n_2, \dots , n_k}$ is given by the following formula
\begin{eqnarray} \label{overlap_dist_min}
\mathbb{P} \left( O_{n_1, n_2, \dots , n_k} > t  \right) &=&  \mathbb{E} \left[  \prod^{k}_{j=1} \overline{G}\left( t + (A_{n_k} - A_{n_j}) \right) \right]  .
\end{eqnarray}
\begin{proof}
\begin{eqnarray*}
\lefteqn{ \mathbb{P} \left( O_{n_1, n_2, \dots , n_k} > t  \right)  } \\
&=& \mathbb{P} \left( ( S_{n_1}  + A_{n_1} - A_{n_k} )^+ \boldsymbol{\wedge} \cdots \boldsymbol{\wedge} ( S_{n_j}  + A_{n_j} - A_{n_k} )^+ \boldsymbol{\wedge} S_{n_k}  > t  \right) \\
&=& \mathbb{P} \left( S_{n_1}+ A_{n_1} - A_{n_k}  > t, \dots , S_{n_j}  + A_{n_j} - A_{n_k}  > t, S_{n_k} > t   \right)  \\
% &=& \mathbb{P} \left( S_{n_1} > t + \mathcal{A}_{(n_k - n_1)}, \dots , S_{n_j}   > t + \mathcal{A}_{(n_k - n_j)}, S_{n_k}  > t   \right)  \\
&=& \mathbb{E} \left[  \prod^{k}_{j=1} \overline{G}\left( t + (A_{n_k} - A_{n_j}) \right)  \right] .
\end{eqnarray*}
This completes the proof. 
\end{proof}
\end{theorem}

\begin{corollary}
Let $O_{n_1, n_2, \dots , n_k}$ be the distribution of the overlap time for customers within distance $n_k - n_1$ in the $GI/M/\infty$ queue, then the tail distribution of $O_{n_1, n_2, \dots , n_k}$ is given by the following formula
\begin{eqnarray} \label{overlap_dist_min}
\mathbb{P} \left( O_{n_1, n_2, \dots , n_k} > t  \right) &=&   e^{- \mu k t  } \mathbb{E} \left[  \prod^{k-1}_{j=1} e^{- \mu j (  A_{n_{j+1}} - A_{n_j}) }  \right].
\end{eqnarray}
\begin{proof}
Using the result from Theorem \ref{tuple_result}, we have 
\begin{eqnarray*}
\mathbb{P} \left( O_{n_1, n_2, \dots , n_k} > t  \right)  
&=& \mathbb{P} \left( S_{n_1}+ A_{n_1} - A_{n_k}  > t, \dots , S_{n_j}  + A_{n_j} - A_{n_k}  > t, S_{n_k} > t   \right)  \\
&=&  e^{- \mu k t  } \mathbb{E} \left[  \prod^{k-1}_{j=1} e^{- \mu (  A_{n_k} - A_{n_j}) }  \right] \\
&=&  e^{- \mu k t  } \mathbb{E} \left[  \prod^{k-1}_{j=1} e^{- \mu j (  A_{n_{j+1}} - A_{n_j}) }  \right] .
\end{eqnarray*}
This completes the proof. 
\end{proof}
\end{corollary}

\begin{corollary}
Let $O_{n_1, n_2, \dots , n_k}$ be the distribution of the overlap time for customers within distance $n_k - n_1$ in the $M/M/\infty$ queue, then the tail distribution of $O_{n_1, n_2, \dots , n_k}$ is given by the following formula
\begin{eqnarray} \label{overlap_dist_min}
\mathbb{P} \left( O_{n_1, n_2, \dots , n_k} > t  \right) &=&  e^{-\mu k t} \prod^{k-1}_{j=1} \left( \frac{\lambda}{\lambda + \mu j } \right)^{n_{j+1} - n_j}  .
\end{eqnarray}
\begin{proof}
Using the result from Theorem \ref{tuple_result}, we have 
\begin{eqnarray*}
\mathbb{P} \left( O_{n_1, n_2, \dots , n_k} > t  \right)   &=&  e^{- \mu k t  } \mathbb{E} \left[  \prod^{k-1}_{j=1} e^{- \mu j (  A_{n_{j+1}} - A_{n_j}) }  \right] \\
&=&  e^{-\mu k t} \prod^{k-1}_{j=1} \left( \frac{\lambda}{\lambda + \mu j } \right)^{n_{j+1} - n_j}  .
\end{eqnarray*}
This completes the proof. 
\end{proof}
\end{corollary}

The first thing to notice is that our results match the two customer case and yield the previous results by \citet{pender2021overlap}.  We also observe that the functional form of the result is not affected by the arrival distribution and is completely governed by the service distribution, which is exponential in this case.  Moreover, the more customers you want to overlap, the smaller the probability since the rate in the exponential is increased according to the number of overlapping customers.  

\subsection{First of Each Batch}
Now that we have analyzed the k-tuple setting in the individual setting i.e. the batch size is equal to 1, we now move to the setting where batches can have an arbitrary size.  We will also perform our analysis for the two batch settings we define in the previous sections i.e. first of each batch and last of each batch.  

\begin{theorem} \label{batch_first_tuple}
Under the first of each batch to leave setting, we let $O_{n_1, n_2, \dots , n_k}$ be the distribution of the overlap time for customers within distance $n_k - n_1$ in the $GI^B/GI/\infty$ queue, then the tail distribution of $O_{n_1, n_2, \dots , n_k}$ is given by the following formula
\begin{eqnarray} \label{overlap_dist_min}
\mathbb{P} \left( O_{n_1, n_2, \dots , n_k} > t  \right) &=&  \mathbb{E} \left[ \prod^{k}_{i=1}  \overline{G} \left( t -  (A_{n_k} - A_{n_i}) \right)^\mathcal{B} \right]     .
\end{eqnarray}
\begin{proof}
\begin{eqnarray*}
\lefteqn{ \mathbb{P} \left( O_{n_1, n_2, \dots , n_k} > t  \right)  } \\
&=& \mathbb{P} \left( ( \min \left( \vec{S}_{n_1}^{B_{n_1}} \right)  + A_{n_1} - A_{n_k} )^+ \boldsymbol{\wedge} \cdots \boldsymbol{\wedge} ( \min \left( \vec{S}_{n_j}^{B_{n_j}} \right)  + A_{n_j} - A_{n_k} )^+ \boldsymbol{\wedge} \min \left( \vec{S}_{n_k}^{B_{n_k}} \right)  > t  \right) \\
&=& \mathbb{P} \left(  \min \left( \vec{S}_{n_1}^{B_{n_1}} \right) + A_{n_1} - A_{n_k}  > t, \dots , \min \left( \vec{S}_{n_j}^{B_{n_j}} \right)  + A_{n_j} - A_{n_k}  > t, \min \left( \vec{S}_{n_k}^{B_{n_k}} \right)  > t   \right)  \\
&=& \mathbb{P} \left( \min \left( \vec{S}_{n_1}^{B_{n_1}} \right)  > t + A_{n_k} - A_{n_1}, \dots , \min \left( \vec{S}_{n_j}^{B_{n_j}} \right)   > t + A_{n_k} - A_{n_j}, \min \left( \vec{S}_{n_k}^{B_{n_k}} \right) > t   \right)  \\
&=& \mathbb{E} \left[ \prod^{k}_{i=1} \left( \sum^{\infty}_{j=1} p_j \cdot \overline{G} \left( t -  (A_{n_k} - A_{n_i}) \right)^j \right) \right] \\
&=& \mathbb{E} \left[ \prod^{k}_{i=1}  \overline{G} \left( t -  (A_{n_k} - A_{n_i}) \right)^\mathcal{B} \right] .
\end{eqnarray*}
This completes the proof. 
\end{proof}
\end{theorem}

\begin{corollary}
Under the first of each batch to leave setting, we let $O_{n_1, n_2, \dots , n_k}$ be the distribution of the overlap time for customers within distance $n_k - n_1$ in the $GI^B/M/\infty$ queue, then the tail distribution of $O_{n_1, n_2, \dots , n_k}$ is given by the following formula
\begin{eqnarray} \label{overlap_dist_min}
\mathbb{P} \left( O_{n_1, n_2, \dots , n_k} > t  \right) &=&   \mathbb{E} \left[ e^{- \mu  k t \mathcal{B}  }   \prod^{k-1}_{i=1} e^{- \mu \mathcal{B} (  A_{n_k} - A_{n_i}) }  \right] .
\end{eqnarray}
\begin{proof}
Using the result from Theorem \ref{batch_first_tuple}, we have that
\begin{eqnarray*}
\mathbb{P} \left( O_{n_1, n_2, \dots , n_k} > t  \right)  
% &=& \mathbb{P} \left( ( \min \left( \vec{S}_{n_1}^{B_{n_1}} \right)  + A_{n_1} - A_{n_k} )^+ \boldsymbol{\wedge} \cdots \boldsymbol{\wedge} ( \min \left( \vec{S}_{n_j}^{B_{n_j}} \right)  + A_{n_j} - A_{n_k} )^+ \boldsymbol{\wedge} \min \left( \vec{S}_{n_k}^{B_{n_k}} \right)  > t  \right) \\
% &=& \mathbb{P} \left(  \min \left( \vec{S}_{n_1}^{B_{n_1}} \right) + A_{n_1} - A_{n_k}  > t, \dots , \min \left( \vec{S}_{n_j}^{B_{n_j}} \right)  + A_{n_j} - A_{n_k}  > t, \min \left( \vec{S}_{n_k}^{B_{n_k}} \right)  > t   \right)  \\
% &=& \mathbb{P} \left( \min \left( \vec{S}_{n_1}^{B_{n_1}} \right)  > t + A_{n_k} - A_{n_1}, \dots , \min \left( \vec{S}_{n_j}^{B_{n_j}} \right)   > t + A_{n_k} - A_{n_j}, \min \left( \vec{S}_{n_k}^{B_{n_k}} \right) > t   \right)  \\
&=& \mathbb{E} \left[  \sum^{\infty}_{j=1} p_j \left( \prod^{k-1}_{i=1} e^{- \mu j ( t +  A_{n_k} - A_{n_i}) }  \cdot e^{- \mu j t  } \right) \right] \\
&=& \sum^{\infty}_{j=1} p_j \cdot e^{- \mu j k t  } \mathbb{E} \left[  \prod^{k-1}_{i=1} e^{- \mu j (  A_{n_k} - A_{n_i}) }  \right] \\
&=& \mathbb{E} \left[ e^{- \mu  k t \mathcal{B}  }   \prod^{k-1}_{i=1} e^{- \mu \mathcal{B} (  A_{n_k} - A_{n_i}) }  \right] .
\end{eqnarray*}
This completes the proof. 
\end{proof}
\end{corollary}

\begin{corollary}
Under the first of each batch to leave setting, we let $O_{n_1, n_2, \dots , n_k}$ be the distribution of the overlap time for customers within distance $n_k - n_1$ in the $GI^B/M/\infty$ queue, then the mean of $O_{n_1, n_2, \dots , n_k}$ is given by the following formula
    \begin{eqnarray} \label{overlap_dist_mean}
\mathbb{E} \left[ O_{n_1, n_2, \dots , n_k} \right] &=&  \frac{1}{\mu k } \mathbb{E} \left[ \mathcal{B}^{-1}  \prod^{k-1}_{i=1} e^{- \mu \mathcal{B} (  A_{n_k} - A_{n_i}) }  \right] .
\end{eqnarray}
\begin{proof}
Using the result from Theorem \ref{batch_first_tuple}, we have that
\begin{eqnarray*}
\mathbb{E} \left[ O_{n_1, n_2, \dots , n_k} \right] &=&  \int^{\infty}_{0} \mathbb{P} \left( O_{n_1, n_2, \dots , n_k} > t  \right)  dt \\
&=& \int^{\infty}_{0}  \mathbb{E} \left[ e^{- \mu  k t \mathcal{B}  }   \prod^{k-1}_{i=1} e^{- \mu \mathcal{B} (  A_{n_k} - A_{n_i}) }  \right] dt \\
&=& \frac{1}{\mu k } \mathbb{E} \left[ \mathcal{B}^{-1}  \prod^{k-1}_{i=1} e^{- \mu \mathcal{B} (  A_{n_k} - A_{n_i}) }  \right] .
\end{eqnarray*}
This completes the proof. 
\end{proof}
\end{corollary}

\subsection{Last of Each Batch}

\begin{theorem} \label{batch_last_tuple}
Under the last of each batch to leave setting, we let $O_{n_1, n_2, \dots , n_k}$ be the distribution of the overlap time for customers within distance $n_k - n_1$ in the $GI^B/GI/\infty$ queue, then the tail distribution of $O_{n_1, n_2, \dots , n_k}$ is given by the following formula
\begin{eqnarray} \label{overlap_dist_max}
\mathbb{P} \left( O_{n_1, n_2, \dots , n_k} > t  \right) &=& 
% \mathbb{E} \left[ \prod^{k}_{i=1} \left( 1 -   \sum^{\infty}_{j=1} p_j \left(  ( 1 - \overline{G} ( t - ( A_{n_k} - A_{n_i} ) )  \right)^j   \right) \right] \\
\mathbb{E} \left[ \prod^{k}_{i=1} \left( 1 -    G \left( t - ( A_{n_k} - A_{n_i} )  \right)^\mathcal{B}   \right) \right] .
\end{eqnarray}
\begin{proof}
\begin{eqnarray*}
\lefteqn{ \mathbb{P} \left( O_{n_1, n_2, \dots , n_k} > t  \right)  } \\
&=& \mathbb{P} \left( ( \max \left( \vec{S}_{n_1}^{B_{n_1}} \right)  + A_{n_1} - A_{n_k} )^+ \boldsymbol{\wedge} \cdots \boldsymbol{\wedge} ( \max \left( \vec{S}_{n_j}^{B_{n_j}} \right)  + A_{n_j} - A_{n_k} )^+ \boldsymbol{\wedge} \max \left( \vec{S}_{n_k}^{B_{n_k}} \right)  > t  \right) \\
&=& \mathbb{P} \left(  \max \left( \vec{S}_{n_1}^{B_{n_1}} \right) + A_{n_1} - A_{n_k}  > t, \dots , \max \left( \vec{S}_{n_j}^{B_{n_j}} \right)  + A_{n_j} - A_{n_k}  > t, \max \left( \vec{S}_{n_k}^{B_{n_k}} \right)  > t   \right)  \\
&=& \mathbb{P} \left( \max \left( \vec{S}_{n_1}^{B_{n_1}} \right)  > t + A_{n_k} - A_{n_1}, \dots , \max \left( \vec{S}_{n_j}^{B_{n_j}} \right)   > t + A_{n_k} - A_{n_j}, \max \left( \vec{S}_{n_k}^{B_{n_k}} \right) > t   \right)  \\
&=& \mathbb{E} \left[ \prod^{k}_{i=1} \left( 1 -   \sum^{\infty}_{j=1} p_j \left(  ( 1 - \overline{G} ( t - ( A_{n_k} - A_{n_i} ) )  \right)^j   \right) \right] \\
&=& \mathbb{E} \left[ \prod^{k}_{i=1} \left( 1 -    G \left( t - ( A_{n_k} - A_{n_i} )  \right)^\mathcal{B}   \right) \right] .
\end{eqnarray*}
This completes the proof. 
\end{proof}
\end{theorem}

\begin{corollary}
Under the last of each batch to leave setting, we let $O_{n_1, n_2, \dots , n_k}$ be the distribution of the overlap time for customers within distance $n_k - n_1$ in the $GI^B/M/\infty$ queue, then the tail distribution of $O_{n_1, n_2, \dots , n_k}$ is given by the following formula
\begin{eqnarray} \label{overlap_dist_max}
\mathbb{P} \left( O_{n_1, n_2, \dots , n_k} > t  \right) &=&  1 - \mathbb{E} \left[   \left( \prod^{k}_{i=1} \left( 1- e^{- \mu \mathcal{B} \left( t +  A_{n_k} - A_{n_i} \right) } \right)   \right) \right]   .
\end{eqnarray}
\begin{proof}
Using the result from Theorem \ref{batch_last_tuple}, we have that
\begin{eqnarray*}
\lefteqn{ \mathbb{P} \left( O_{n_1, n_2, \dots , n_k} > t  \right)  } \\
&=& \mathbb{P} \left( ( \max \left( \vec{S}_{n_1}^{B_{n_1}} \right)  + A_{n_1} - A_{n_k} )^+ \boldsymbol{\wedge} \cdots \boldsymbol{\wedge} ( \max \left( \vec{S}_{n_j}^{B_{n_j}} \right)  + A_{n_j} - A_{n_k} )^+ \boldsymbol{\wedge} \max \left( \vec{S}_{n_k}^{B_{n_k}} \right)  > t  \right) \\
&=& \mathbb{P} \left(  \max \left( \vec{S}_{n_1}^{B_{n_1}} \right) + A_{n_1} - A_{n_k}  > t, \dots , \max \left( \vec{S}_{n_j}^{B_{n_j}} \right)  + A_{n_j} - A_{n_k}  > t, \max \left( \vec{S}_{n_k}^{B_{n_k}} \right)  > t   \right)  \\
&=& \mathbb{P} \left( \max \left( \vec{S}_{n_1}^{B_{n_1}} \right)  > t + A_{n_k} - A_{n_1}, \dots , \max \left( \vec{S}_{n_j}^{B_{n_j}} \right)   > t + A_{n_k} - A_{n_j}, \max \left( \vec{S}_{n_k}^{B_{n_k}} \right) > t   \right)  \\
&=& 1 - \mathbb{E} \left[  \sum^{\infty}_{j=1} p_j \left( \prod^{k}_{i=1} \left( 1- e^{- \mu j \left( t +  A_{n_k} - A_{n_i} \right) } \right)   \right) \right] \\
&=& 1 - \mathbb{E} \left[   \left( \prod^{k}_{i=1} \left( 1- e^{- \mu \mathcal{B} \left( t +  A_{n_k} - A_{n_i} \right) } \right)   \right) \right] .
\end{eqnarray*}
This completes the proof. 
\end{proof}
\end{corollary}

\section{Conclusion} \label{conclusion}

In this paper, we consider the overlap times for customers in the $GI^B/GI/\infty$ queue.   We provide an analysis from an individual and a batch perspective.  Despite our analysis, there are many avenues for additional research.  First, as we mentioned before, there are several ways to analyze the batch perspective.  One such way would be to explicitly analyze the order statistics of each batch.  In particular one might be interested in when a certain percentage of customers have departed from each batch.  

Second, one can also extend our work by considering generalizations such as dependent arrivals or service times like in \citet{pang2012impact, pang2013two, daw2018queues}.  In particular, dependent service times would impact the tail distribution quite significantly since it no longer splits into a product of individual distributions.  It would also be interesting to consider the situation where there are a finite number of servers or abandonments to assess the impact on the tail distribution of the overlap time, see for example \citet{ko2022overlapping}.  Finally, we are also interested in designing queueing systems to minimize overlap times under some cost criterion.  We hope to consider these important extensions in future work.  

\section*{Acknowledgements}
Jamol Pender would like to acknowledge the gracious support of the National Science Foundation DMS Award \# 2206286.  Sergio Palomo would like to acknowledge the gracious support of the Sloan Foundation for supporting his graduate studies.  

\bibliographystyle{plainnat}
\bibliography{references}
\end{document}